# Maxiset in sup-norm for kernel estimators

KARINE BERTIN

*Departamento de Estadística, Universidad de Valparaíso, DIPUV 33/05, FONDECYT 1061184*

VINCENT RIVOIRARD

*Laboratoire de Mathématiques, CNRS UMR 8628, Université Paris Sud*

**ABSTRACT.** In the Gaussian white noise model, we study the estimation of an unknown multidimensional function $f$ in the uniform norm by using kernel methods. The performances of procedures are measured by using the maxiset point of view: we determine the set of functions which are well estimated (at a prescribed rate) by each procedure. So, in this paper, we determine the maxisets associated to kernel estimators and to the Lepski procedure for the rate of convergence of the form $(\log n/n)^{\beta/(2\beta+d)}$. We characterize the maxisets in terms of Besov and Hölder spaces of regularity $\beta$.

*Key Words:* Besov spaces, kernel estimator, Lepski procedure, maxiset, Uniform norm.

## 1 Introduction

We consider the maxiset point of view in the classical Gaussian white noise model

$$dY_t = f(t)dt + \frac{\sigma}{\sqrt{n}}dW_t, \quad t \in [0,1]^d, \tag{1}$$

where $f : \mathbb{R}^d \to \mathbb{R}$ is an unknown function, $W$ is the Brownian sheet in $[0,1]^d$, $\sigma > 0$ is known and $n \in \mathbb{N}^*$. We take a level of noise of the form $\sigma/\sqrt{n}$ to refer to the equivalence between the model (1) and the regression model with $n$ observations and noise level equal to $\sigma$. In this paper, we study the estimation of $f$ on $[0,1]^d$ from the observations $\{Y_t, t \in [0,1]^d\}$. For this purpose, we assume that $f$ belongs to $\mathbb{L}_\infty^{per}(\mathbb{R}^d)$ the set of 1-periodic functions that belong to $\mathbb{L}_\infty(\mathbb{R}^d)$. The quality of an estimator $\hat{f}_n$ is characterized by its risk in sup-norm

$$R_n(\hat{f}_n) = \mathbb{E}\left(\|\hat{f}_n - f\|_\infty^p\right),$$

where $\|g\|_\infty = \operatorname{ess\,sup}_{t \in [0,1]^d} |g(t)|$ and $p \geq 1$.

In a general way, for non parametric framework, there are three steps for the statistician when he faces the problem of estimating $f$: the choice of the method (kernel, Fourier series, wavelet,...), the determination of parameters of the method (the bandwidth, the number of coefficients that have to be estimated,...) and the evaluation of the quality of his procedure $\hat{f} = (\hat{f}_n)_n$ (the word "procedure" sets the couple method-parameters) by computing the rate of convergence of $R_n(\hat{f}_n)$. In the non-parametric setting, the minimax theory is the most popular point of view: it consists in choosing a functional space $\mathcal{F} \subset \mathbb{L}_\infty^{per}(\mathbb{R}^d)$ and ensuring that $\hat{f}$ achieves the best rate on $\mathcal{F}$. But, at first, the rate could be unknown, secondly, the choice of $\mathcal{F}$ is arbitrary (what kind of spaces has to be considered: Sobolev spaces? Besov spaces? why?), thirdly, $\mathcal{F}$ could contain very bad functions $g$ (in the sense that $g$ is difficult to estimate). Since the unknown quantity $f$ could be easier to estimate, the used procedure could be too pessimistic and not adapted to the data. More embarrassing in practice, several minimax procedures may be proposed and the practitioner has no way to decide if he has no practical experience. To answer these practical questions, an other point of view has recently appeared: the maxiset point of view (see for instance Kerkyacharian & Picard (2000)). It consists in deciding the accuracy of the estimate by fixing a prescribed rate $\rho$ and to point out all the functions $f$ such that $f$ can be estimated by the procedure $\hat{f}$ at the target rate $\rho$. The maxiset of the procedure $\hat{f}$ for this rate $\rho$ is the set of all these functions. So, in our framework, we set the following definition.

**Definition 1.** *Let $1 \leq p < \infty$, $\rho = (\rho_n)_n$ a decreasing sequence of positive real numbers and let $\hat{f} = (\hat{f}_n)_n$ be an estimation procedure. The maxiset of $\hat{f}$ associated with the rate $\rho$ and the $\|.\|_\infty^p$-loss is:*

$$MS(\hat{f}, \rho, p) = \left\{ f \in \mathbb{L}_\infty^{per}(\mathbb{R}^d) : \quad \sup_n \left[ \rho_n^{-p} \mathbb{E} \left( \|\hat{f}_n - f\|_\infty^p \right) \right] < \infty \right\}.$$

The maxiset point of view brings answers to the previous questions. Indeed, there is no a priori functional assumption and then, the practitioner does not need to restrict his study to an arbitrary functional space. The practitioner states the desired accuracy and then, knows the quality of the used procedure. Obviously, he chooses the procedure with the largest maxiset since the larger the maxiset the better the procedure. Previous results concerning the maxiset approach are the following. Maxisets of linear procedures are Besov spaces $\mathcal{B}_{q,\infty}^\beta$ when investigated under the $\mathbb{L}_q$-norm ($1 < q < \infty$) and with polynomial rates of convergence (see Kerkyacharian & Picard (1993)). These results have been generalized by Rivoirard (2004) who proved that linear procedures are outperformed by non linear ones from the maxiset approach. Maxisets for adaptive local and global thresholding rules and Bayesian rules have been investigated in different settings (in the white noise model, in density estimation, in inverse problems or for integrated quadratic functionals estimation). See Kerkyacharian & Picard (2000), Cohen *et al* (2001), Kerkyacharian & Picard (2002), Rivoirard (2005), Autin (2006), Autin *et al* (2006) and Rivoirard & Tribouley (2006). All these results based on wavelet procedures have been derived for the $\mathbb{L}_2$-norm or for the $\mathbb{L}_q$-norm ($1 < q < \infty$). Furthermore, most of these maxiset results are also established for one-dimensional functions and with the rate $\psi(\beta) = (\psi_n(\beta))_n$ where

$$\psi_n(\beta) = \left(\frac{\log n}{n}\right)^{\beta/(2\beta+d)}$$

and the maxisets are not Besov spaces $\mathcal{B}_{q,\infty}^\beta$ but Lorentz spaces that can be viewed as weak versions of Besov spaces and are strictly larger than Besov spaces. So, the framework of this paper is quite different since, for estimating multidimensional functions, we consider kernel estimators and the $\mathbb{L}_\infty$-norm.

In the non-adaptive minimax framework, estimation in sup-norm has been studied by Ibragimov & Hasminskii (1981) for one-dimensional Hölder functions and by Stone (1982) for multidimensional isotropic Hölder functions. They proved that the minimax rate of convergence is $\psi(\beta)$ for estimation of functions with known regularity $\beta$. In the adaptive case, Lepski (1992) and Bertin (2005) obtained the same rate for estimation of Hölderian functions for the one- and multidimensional case. Most of these results are based on kernel rules. Considering Hölderian functions is a very classical choice in this minimax setting. Results using kernels for estimation in sup-norm are reviewed in Tsybakov (2004).

As explained before, our goal in this paper is to investigate maxisets in sup-norm for kernel procedures and we consider the rate $\psi(\beta)$. Note that it is of interest in using the sup-norm in estimation since it provides a band of confidence around the estimator. When the $\mathbb{L}_q$-norm is considered with $q < \infty$, maxisets of classical wavelet estimators for the rate $\psi(\beta)$ are not classical Besov spaces (see previously). So the question is the following. Are the maxisets obtained in the framework of this paper classical Besov spaces? Roughly speaking, the answer is yes, as shown by the following precise description of our results. To derive maxisets for classical kernel estimators, we have to prove two inclusions. In Theorem 1, we first prove that the maxiset of a very general class of kernel estimators defined in (2) is included in the set of functions that satisfy Condition (4). Then, in Theorem 2, under this condition, the converse inclusion is proved for a more specific class of kernel estimators defined in (6). Theorem 3 gathers the two previous results giving exactly what is the maxiset of this specific class of estimators. Actually, using the smoothing operator $f \mapsto K_{h_n} * f$ where $K_{h_n}$ is the kernel from which our procedure is defined, the maxiset is the set of the functions $f$ that can be approximated by $K_{h_n} * f$ at the rate $\psi_n(\beta)$ (see (10) for a precise definition of this set). In Theorem 4, under some conditions on the kernel, we prove that this set is in fact the Besov space $\mathcal{B}^{\beta}_{\infty,\infty}$ when $\beta$ is not an integer. Since, in this case, $\mathcal{B}^{\beta}_{\infty,\infty}$ is equal to the set of $\beta$-Hölder functions, this result justifies the classical choice of Hölder spaces to study minimax properties of procedures under the $\|\cdot\|_\infty$-loss. When $\beta$ is an integer, the maxiset contains the $\beta$-Hölder functions and is included in $\mathcal{B}^{\beta}_{\infty,\infty}$. In fact, it was already known that $\beta$-Hölder functions can be estimated at the rate $\psi_n(\beta)$ (see previously), but, roughly speaking, we prove that these functions are the only ones and this result is new.

For the previous results, our kernel procedures depend on $\beta$ (see (6)). In the minimax approach, as said previously, an adaptive estimation procedure has been proposed by Lepski (1992) that achieves the same minimax properties as non-adaptive kernel procedures when Hölder functions are considered. Naturally, our next goal is to compare the maxiset performance of the Lepski procedure with previous kernel procedures from the maxiset point of view. In Theorem 5, we prove that, for all $\beta \notin \mathbb{N}$, the maxiset associated to the Lepski procedure is the Besov space $\mathcal{B}^{\beta}_{\infty,\infty}$. So, the adaptive Lepski procedure achieves exactly the same maxiset performance as estimators defined in (6).

The paper is organized as follows. Section 2 contains results described previously. More precisely, in Subsection 2.1, we introduce kernel estimators. In Subsection 2.2, maxisets of these estimators are derived. Their characterizations in terms of functional spaces are given in Subsection 2.3. The Lepski procedure is introduced in Subsection 2.4 and its maxisets are pointed out. Section 3 is devoted to the proofs.

## 2 Main results

### 2.1 Kernel estimators

We consider the classical class $\mathcal{I}$ of kernel estimators.

**Definition 2.** *The class $\mathcal{I}$ is the class of kernel estimators $\left(\tilde{f}_{K,h_n}\right)_{n \in \mathbb{N}^*}$ of the form:*

$$\tilde{f}_{K,h_n}(t) = \frac{1}{h_n^d} \int_{\mathbb{R}^d} K\left(\frac{t-u}{h_n}\right) dY_u, \quad t \in [0,1]^d \tag{2}$$

*where $(h_n)_{n \in \mathbb{N}^*}$ is a sequence of positive numbers that tends to 0 and $K$ is a function $K : \mathbb{R}^d \to \mathbb{R}$ that satisfies the two following conditions:*

$(A_1)$ *$K$ has a compact support,*

$(A_2)$ *$\|K\|_2^2 = \int_{\mathbb{R}^d} K^2(u) du < \infty$.*

Note that since we consider some $K$ compactly supported and that in the sequel, the so-called bandwidth parameter $h_n$ will be small, the estimator $\tilde{f}_{K,h_n}(t)$ is well defined when $t$ is far from the boundary of $[0,1]^d$ (because $K\left(\frac{t-u}{h_n}\right) = 0$ when $u$ is outside of $[0,1]^d$). However, some border problems arise when we want to define $\tilde{f}_{K,h_n}(t)$ for $t$ close to the boundary of $[0,1]^d$. To answer this issue, for any $j = (j_1, \ldots, j_d) \in \mathbb{Z}^d$ and any function $g \in \mathbb{L}_2(\mathbb{R}^d)$, if $\mathcal{J} = [j_1, j_1+1] \times \cdots \times [j_d, j_d+1]$, we set

$$\int_{\mathcal{J}} g(u) dW_u = \int_{[0,1]^d} g(u+j) dW_u. \tag{3}$$

Since $f$ is 1-periodic, this implies that

$$\int_{\mathcal{J}} g(u) dY_u = \int_{[0,1]^d} g(u+j) dY_u.$$

This trick allows to define integrals of the form $\int_{\mathbb{R}^d} g(u) dY_u$ by decomposing $\mathbb{R}^d$ as a union of compact intervals of the form $[j_1, j_1+1] \times \cdots \times [j_d, j_d+1]$. Using again the periodicity of $f$, we obtain the classical form for $\mathbb{E}\left[\tilde{f}_{K,h_n}(t)\right]$:

$$\forall\, t \in [0,1]^d, \quad \mathbb{E}\left[\tilde{f}_{K,h_n}(t)\right] = \frac{1}{h_n^d} \int_{\mathbb{R}^d} K\left(\frac{t-u}{h_n}\right) f(u) du = K_{h_n} * f(t)$$

where for any $t \in \mathbb{R}^d$,

$$K_{h_n}(t) = \frac{1}{h_n^d} K\left(\frac{t}{h_n}\right)$$

and for any functions $f_1$ and $f_2$, $f_1 * f_2$ denotes the standard convolution product on $\mathbb{R}^d$ between $f_1$ and $f_2$.

Note that these border effects can also be dealt with folded kernels. See for instance, Bertin (2004).

Most of the functions we consider are 1-periodic (in particular since $f$ is 1-periodic, $K_{h_n} * f$ is also 1-periodic) and in this case the sup-norm on $[0,1]^d$ is identical to the sup-norm on $\mathbb{R}^d$.

## 2.2 Maxiset of kernel estimators

Before deriving maxisets for classical kernel estimators, let us point out the following theorem.

**Theorem 1.** *Let $\beta > 0$ and let us consider $(\tilde{f}_{K,h_n})_n \in \mathcal{I}$. If the function $f$ satisfies*

$$\sup_n \left[\psi_n^{-p}(\beta) \times \mathbb{E}\|\tilde{f}_{K,h_n} - f\|_\infty^p\right] < \infty,$$

*then we have*

$$\sup_n \left[h_n^{-\beta} \|\mathbb{E}\tilde{f}_{K,h_n} - f\|_\infty\right] < \infty \tag{4}$$

*and there exists a positive constant $C$ such that $\forall\, n \in \mathbb{N}^*$,*

$$h_n^{-1} \leq C \left(\frac{\log(n)}{n}\right)^{-1/(2\beta+d)}. \tag{5}$$

This result shows that for the rate $\psi(\beta)$, the maxiset of the very classical family of estimators $\mathcal{I}$ cannot be larger than the set of functions that satisfies (4). Furthermore, the bandwidth parameter $h_n$ cannot go to 0 too quickly. Theorem 1 is proved in Section 3 but let us give the main tools that allow to prove this result.

In the minimax setting for an estimator $f_n^*$ and $p \geq 1$, we usually use the classical decomposition of the risk in bias and variance terms:

$$\mathbb{E}\|f_n^* - f\|_\infty^p \leq 2^{p-1}\left(\|\mathbb{E}f_n^* - f\|_\infty^p + \mathbb{E}\|f_n^* - \mathbb{E}f_n^*\|_\infty^p\right).$$

In the maxiset setting, and in particular to prove Theorem 1, we use the following converse result that shows that controlling the risk allows to control the bias and the variance terms.

**Lemma 1.** *For any estimator $f_n^*$, we have:*

$$\|\mathbb{E} f_n^* - f\|_\infty^p \leq \mathbb{E}\|f_n^* - f\|_\infty^p,$$

$$\mathbb{E}\|f_n^* - \mathbb{E} f_n^*\|_\infty^p \leq 2^p \mathbb{E}\|f_n^* - f\|_\infty^p.$$

The proof of Theorem 1 also relies on the following proposition concerning the variance term that actually provides the lower bound for the bandwidth parameter.

**Proposition 1.** *Let us consider $(\tilde{f}_{K,h_n})_n \in \mathcal{I}$. For any $\delta > 0$, there exists $n_0 \in \mathbb{N}$ such that for any $n \geq n_0$,*

$$\mathbb{E}\|\tilde{f}_{K,h_n} - \mathbb{E}\tilde{f}_{K,h_n}\|_\infty^p \geq (1-\delta) \left(\frac{2d\sigma^2 \|K\|_2^2 |\log(h_n)|}{n h_n^d}\right)^{p/2}.$$

Now, our goal is to build a procedure achieving the rate $\psi(\beta)$ under the constraints (4) and (5). Such a procedure is built in Theorem 2. In particular (7) provides the optimal choice for the bandwidth parameter.

**Theorem 2.** *Let $\beta > 0$. We consider the estimator $\hat{f}_{n,\beta}$ defined for $t \in [0,1]^d$ by*

$$\hat{f}_{n,\beta}(t) = \frac{1}{h_{n,\beta}^d} \int_{\mathbb{R}^d} K\left(\frac{t-u}{h_{n,\beta}}\right) dY_u, \tag{6}$$

*with*

$$h_{n,\beta} = C \left(\frac{\log n}{n}\right)^{\frac{1}{2\beta+d}} \tag{7}$$

*where $C$ is a positive constant. We suppose that $K$ satisfies $(A_1)$, $(A_2)$ and the following condition:*

$(A_3)$ *for all $t \in \mathbb{R}^d$ such that $\|t\| \leq 1$, $\int_{\mathbb{R}^d}(K(t+u) - K(u))^2 du \leq C\|t\|^{2\gamma}$, where $\|\cdot\|$ is a norm of $\mathbb{R}^d$, $C$ is a positive constant and $\gamma \in (0,1]$.*

*Then if $f$ satisfies*

$$\sup_n \left[h_{n,\beta}^{-\beta} \|\mathbb{E}\hat{f}_{n,\beta} - f\|_\infty\right] < \infty, \tag{8}$$

*we have*

$$\sup_n \left[ \psi_n^{-p}(\beta) \mathbb{E} \left[ \|\hat{f}_{n,\beta} - f\|_\infty^p \right] \right] < \infty. \tag{9}$$

Now, Theorems 1 and 2 easily imply the following maxiset result.

**Theorem 3.** *For any $\beta > 0$, let us set*

$$\hat{f}_\beta = (\hat{f}_{n,\beta})_n,$$

*where for any $n$, $\hat{f}_{n,\beta}$ is defined by (6) with the bandwidth parameter defined in (7) and such that the kernel $K$ satisfies $(A_1)$, $(A_2)$ and $(A_3)$. Then,*

$$MS(\hat{f}_\beta, \psi(\beta), p) = \left\{ f \in \mathbb{L}_\infty^{per}(\mathbb{R}^d) : \sup_n \left[ h_{n,\beta}^{-\beta} \|K_{h_{n,\beta}} * f - f\|_\infty \right] < \infty \right\}. \tag{10}$$

Note that $MS(\hat{f}_\beta, \psi(\beta), p)$ does not depend on the parameter $p$. But this maxiset depends on the kernel $K$ and on the bandwidth parameter $h_{n,\beta}$. Furthermore, $MS(\hat{f}_\beta, \psi(\beta), p)$ does not look like a classical functional space. In the following subsection, by adding some mild conditions, we characterize this maxiset by classical functional spaces.

## 2.3  Characterization of the maxiset in terms of functional spaces

### 2.3.1  Functional classes

Let us recall the definition of some classical functional spaces that will play a capital role in the maxiset setting of this paper. First, let us adopt Meyer's point of view to introduce Besov spaces viewed as approximation spaces (see Meyer (1990) p. 49). This approach is natural in the context of this paper.

**Definition 3.** *For $\beta > 0$, the Besov space $\mathcal{B}_{\infty,\infty}^{per,\beta}$ is the set of functions $f$ that belong to $\mathbb{L}_\infty^{per}(\mathbb{R}^d)$ and satisfy the following property: for any integer $N > \beta$, there exists a sequence of functions $(f_j)_{j \in \mathbb{N}^*}$ belonging to $\mathbb{L}_\infty^{per}(\mathbb{R}^d)$ such that*

$$\sup_{j \in \mathbb{N}^*} 2^{j\beta} \|f - f_j\|_\infty < \infty$$

and

$$\sup_{j \in \mathbb{N}^*} 2^{-(N-\beta)j} \left\| \frac{\partial^{\alpha_1}}{\partial t_1^{\alpha_1}} \cdots \frac{\partial^{\alpha_d}}{\partial t_d^{\alpha_d}} f_j \right\|_\infty < \infty,$$

where $\alpha = (\alpha_1, \ldots, \alpha_d) \in \mathbb{N}^d$ satisfies $\sum_{i=1}^d \alpha_i = N$.

Now, let us introduce Hölder spaces.

**Definition 4.** *For $\beta \in (0, 1]$, the Hölder space $\Sigma^{per}(\beta)$ is the set of continuous functions $f$ that belong to $\mathbb{L}_\infty^{per}(\mathbb{R}^d)$ and satisfy:*

$$\sup_{x \neq y} \frac{|f(x) - f(y)|}{\|x - y\|^\beta} < \infty.$$

*For $\beta > 1$, the Hölder space $\Sigma^{per}(\beta)$ is the set of functions $f$ of class $C^m$ that belong to $\mathbb{L}_\infty^{per}(\mathbb{R}^d)$ and such that all the derivatives of order $m$ belong to $\Sigma^{per}(\alpha)$ where $m = \lfloor \beta \rfloor = \max\{l \in \mathbb{N}, l < \beta\}$ and $\beta = m + \alpha$.*

When $\beta \notin \mathbb{N}^*$ the Hölder space $\Sigma^{per}(\beta)$ and the Besov space $\mathcal{B}_{\infty,\infty}^{per,\beta}$ are identical (see for instance Meyer (1990) p. 52–53). This is not true when $\beta$ is an integer and $\Sigma^{per}(\beta)$ is strictly included in $\mathcal{B}_{\infty,\infty}^{per,\beta}$. We have the following result proved by Meyer (1990) (cf. p. 53).

**Proposition 2.** *For $\beta \in (0, 1]$, a continuous function $f \in \mathbb{L}_\infty^{per}(\mathbb{R}^d)$ belongs to the Besov space $\mathcal{B}_{\infty,\infty}^{per,\beta}$ if and only if*

- *when $0 < \beta < 1$,*

$$\sup_{x \neq y} \frac{|f(x) - f(y)|}{\|x - y\|^\beta} < \infty,$$

- *when $\beta = 1$,*

$$\sup_{x \in \mathbb{R}^d, y \neq 0} \frac{|f(x + y) + f(x - y) - 2f(x)|}{\|y\|} < \infty.$$

*For $\beta > 1$, the Besov space $\mathcal{B}_{\infty,\infty}^{per,\beta}$ is the set of functions $f$ of class $C^m$ that belong to $\mathbb{L}_\infty^{per}(\mathbb{R}^d)$ and such that all the derivatives of order $m$ belong to $\mathcal{B}_{\infty,\infty}^{per,\alpha}$ where $m = \lfloor \beta \rfloor = \max\{l \in \mathbb{N}, l < \beta\}$ and $\beta = m + \alpha$.*

In the following sections, to avoid tedious notations, we note $\mathcal{B}^{\beta}_{\infty,\infty}$ instead of $\mathcal{B}^{per,\beta}_{\infty,\infty}$ and $\Sigma(\beta)$ instead of $\Sigma^{per}(\beta)$. Besov spaces and Hölder spaces will naturally characterize maxisets of the kernel procedures.

### 2.3.2 Hypothesis on the kernel estimators

Before characterizing maxisets for kernel estimators, we need to restrict the class $\mathcal{I}$. For this purpose, let us introduce $\mathcal{I}(N)$ defined as follows.

**Definition 5.** *For $N \in \mathbb{N}^*$, $\mathcal{K}(N)$ is the set of the functions $K : \mathbb{R}^d \to \mathbb{R}$ that satisfy conditions $(A_1)$, $(A_2)$, $(A_3)$ and*

$(A_4)$ $\int_{\mathbb{R}^d} K(u) du = 1$,

$(A_5)$ *for any $(\alpha_1, \ldots, \alpha_d) \in \mathbb{N}^d$, such that $\sum_{i=1}^d \alpha_i \leq N$, we have*

$$\int_{\mathbb{R}^d} \left| \frac{\partial^{\alpha_1}}{\partial t_1^{\alpha_1}} \cdots \frac{\partial^{\alpha_d}}{\partial t_d^{\alpha_d}} K(t) \right| dt < \infty,$$

$(A_6)$ *for all polynom $P$ of degree inferior to $N$ such that $P(0) = 0$,*

$$\int_{\mathbb{R}^d} P(u) K(u) du = 0.$$

*The set $\mathcal{H}$ is the set of sequences $(h_n)_{n \in \mathbb{N}^*}$ of the form $h_n = 2^{-m_n}$, $n \in \mathbb{N}^*$, with a sequence $(m_n)_{n \in \mathbb{N}^*}$ that satisfies*

1. *$(m_n)_{n \in \mathbb{N}^*}$ is non decreasing,*
2. *$\lim_{n \to +\infty} m_n = +\infty$,*
3. *$\sup_n (m_{n+1} - m_n) < \infty$.*

*For $N \in \mathbb{N}^*$, $\mathcal{I}(N)$ is the class of kernel estimators*

$$\mathcal{I}(N) = \left\{ \left( \tilde{f}_{K,h_n} \right)_{n \in \mathbb{N}^*} : \tilde{f}_{K,h_n}(t) = \frac{1}{h_n^d} \int_{\mathbb{R}^d} K\left( \frac{t-u}{h_n} \right) dY_u, \ t \in [0,1]^d, K \in \mathcal{K}(N), (h_n)_{n \in \mathbb{N}^*} \in \mathcal{H} \right\}.$$

Conditions for belonging to $\mathcal{H}$ are very mild and are for instance satisfied by the sequence $(h_{n,\beta})_n$ introduced in Theorem 2. Furthermore, the sets $\mathcal{K}(N)$ contain most of the kernels used in estimation:

- $K(x) = 1_{[-1/2,1/2]}(x)$ for $N = 1$.

- $K(x) = c_\beta (1 - \sum_{i=1}^{d} |x_i|^\beta)_+$ with $\beta \geq 1$ for $N = 1$, with $\beta \geq 2$ for $N = 2$ where for any $x \in \mathbb{R}$, $x_+ = \max(0, x)$.

- $K(x) = d_\beta (1 - \sum_{i=1}^{d} |x_i|^\beta)_+^2$ with $\beta \geq 1$ for $N = 1$.

- For $N > 2$, see the construction of higher order kernels by Tsybakov (2004) §1.2.2 for instance.

So, the class $\mathcal{I}(N)$ is a very general class of kernel estimators. Note that in Condition $(A_5)$ the kernel $K$ only needs to be differentiable almost everywhere.

### 2.3.3 Characterizations of maxisets for kernel estimators

Using the class $\mathcal{I}(N)$, we prove the following result. For $\beta > 0$, we note

$$\lceil \beta \rceil = \min \{l \in \mathbb{N} : \quad l > \beta\}.$$

**Theorem 4.** *Consider the procedure $\hat{f}_\beta = (\hat{f}_{n,\beta})_n$ defined in Theorem 2 with $K \in \mathcal{K}(\lceil \beta \rceil)$ and $h_{n,\beta} = C \left( \frac{\log n}{n} \right)^{\frac{1}{2\beta + d}}$.*

1. *If $\beta$ is not an integer*

$$MS(\hat{f}_\beta, \psi(\beta), p) = \mathcal{B}_{\infty,\infty}^\beta,$$

2. *if $\beta$ is an integer*

$$\Sigma(\beta) \subset MS(\hat{f}_\beta, \psi(\beta), p) \subset \mathcal{B}_{\infty,\infty}^\beta.$$

Theorem 4 is proved in Section 3 as a consequence of Theorem 3. This result establishes that the set of functions that can be estimated at the classical rate $\psi(\beta)$ is exactly the functions that belong to $\mathcal{B}_{\infty,\infty}^\beta$ when $\beta$ is not an integer. When $\beta$ is an integer, there is a slight ambiguity resulting from the strict inclusion of $\Sigma(\beta)$ in $\mathcal{B}_{\infty,\infty}^\beta$.

Until now, we have investigated maxisets for kernel procedures depending on $\beta$ through the bandwidth parameter. Now, in view of adaptation, the question is the following. Can we build a kernel procedure $\hat{f}$ such that for any $\beta > 0$, $\hat{f}$ achieves the same maxiset properties as the procedure $\hat{f}_\beta$ built in Theorem 2?

This problem is answered in the next section by considering the Lepski procedure. And we prove as previously that, roughly speaking, Besov spaces are maxisets of adaptive kernel procedures.

## 2.4 Maxisets for the Lepski procedure

In this subsection, we determine the maxiset associated to Lepski procedure (Lepski (1992)). Let $B = \{\beta_1, \ldots, \beta_L\}$ a finite subset of $(0, +\infty)^d$ such that $\beta_i < \beta_j$ if $i < j$ and the $\beta_i$'s are non-integer. For each $\beta \in B$, we consider the procedure $\hat{f}_\beta = (\hat{f}_{n,\beta})_n$ defined in Theorem 2 with $K \in \mathcal{K}(\lceil \beta \rceil)$ and $h_{n,\beta} = C \left(\frac{\log n}{n}\right)^{\frac{1}{2\beta+d}}$. We set

$$\hat{\beta} = \max\left\{ u \in B : \quad \|\hat{f}_{n,u} - \hat{f}_{n,\gamma}\|_\infty \leq \eta_n(\gamma), \ \forall \gamma \leq u \right\},$$

with

$$\eta_n(\gamma) = C_1 \psi_n(\gamma),$$

and $C_1$ is a constant assumed to be large enough (cf. Lepski (1992) and Bertin (2005) for a precise choice of the constant $C_1$). Denote this procedure $\hat{f} = (\hat{f}_{n,\hat{\beta}})_n$. The Lepski procedure is based on the fact that while $\gamma \leq \beta \leq \delta$ and $f$ is of regularity $\delta$, the bias of $\hat{f}_{n,\beta} - \hat{f}_{n,\gamma}$ is bounded from above by a term of order $\psi_n(\gamma)$. We have the following theorem.

**Theorem 5.** *Let $\beta \in B$. We have*

$$MS(\hat{f}, \psi(\beta), p) = \mathcal{B}^\beta_{\infty,\infty}.$$

This result proves that the adaptive kernel procedure $\hat{f}$ achieves the same performance as $\hat{f}_\beta$ from the maxiset point of view. But $\hat{f}$ does not depend on $\beta$ and automatically adapts to the unknown regularity of the function to be estimated. To prove Theorem 5, we first use arguments of Bertin (2005) to derive the inclusion $\mathcal{B}^\beta_{\infty,\infty} \subset MS(\hat{f}, \psi(\beta), p)$. The inclusion $MS(\hat{f}, \psi(\beta), p) \subset \mathcal{B}^\beta_{\infty,\infty}$ is expected since we guess that the maxiset performances of $\hat{f}$ cannot be stronger than those of $\hat{f}_\beta$. Technical details of this proof are given in Section 3.

# 3 Proofs

## 3.1 Proof of Lemma 1

The first inequality is a simple consequence of Jensen inequality. Now,

$$\begin{aligned}\mathbb{E}\|f_n^* - \mathbb{E}f_n^*\|_\infty^p &\leq 2^{p-1}\left(\|\mathbb{E}f_n^* - f\|_\infty^p + \mathbb{E}\|f_n^* - f\|_\infty^p\right) \\ &\leq 2^p \mathbb{E}\|f_n^* - f\|_\infty^p,\end{aligned}$$

that gives the result.

## 3.2 Proof of Proposition 1

We set

$$\forall\, t \in [0,1]^d, \quad Z_n(t) = \tilde{f}_{K,h_n}(t) - \mathbb{E}(\tilde{f}_{K,h_n}(t)) = \frac{\sigma}{\sqrt{n h_n^d}} \int_{\mathbb{R}^d} K\left(\frac{t-u}{h_n}\right) dW_u$$

and

$$\forall\, t \in [0,1]^d, \quad \xi_t = \frac{\sqrt{n h_n^d}}{\sigma} Z_n(t) = \frac{1}{h_n^{d/2}} \int_{\mathbb{R}^d} K\left(\frac{t-u}{h_n}\right) dW_u.$$

Let $A > 0$ such that the support of $K$ is included in $[-A, A]^d$. We set,

$$m = \left\lfloor \frac{1}{2Ah_n} \right\rfloor - 1,$$

and

$$\forall\, i = (i_1, i_2, \ldots, i_d) \in \{1, \ldots, m\}^d, \quad t_i = (2i_1 A h_n, 2i_2 A h_n, \ldots, 2i_d A h_n).$$

Since the support of $K$ is included in $[-A, A]^d$, then $(\xi_{t_i})_{i \in \{1,\ldots,m\}^d}$ are independent centered Gaussian variables with common variance $s^2 = \|K\|_2^2$. We also have for any $r > 0$,

$$
\begin{aligned}
\mathbb{P}\left(\sup_{t \in [0,1]^d} |\xi_t| > r\right) &\geq \mathbb{P}\left(\sup_{i \in \{1,\ldots,m\}^d} \xi_{t_i} > r\right) \\
&= 1 - \mathbb{P}\left(\sup_{i \in \{1,\ldots,m\}^d} \xi_{t_i} \leq r\right) \\
&= 1 - \mathbb{P}\left(\xi_{t_1} \leq r\right)^{m^d} \\
&= 1 - (1 - \phi(r/s))^{m^d} \\
&= 1 - \exp(m^d \log(1 - \phi(r/s))) \\
&\geq 1 - \exp\left(-m^d \phi(r/s)\right)
\end{aligned}
$$

where for any $x \in \mathbb{R}_+$,

$$
\phi(x) = \frac{1}{\sqrt{2\pi}} \int_x^{+\infty} \exp\left(-\frac{v^2}{2}\right) dv \geq \frac{1}{\sqrt{2\pi}} \exp\left(-\frac{x^2}{2}\right) \frac{x}{1 + x^2}.
$$

So, there exists $h_0$ such that if $h_n \leq h_0$,

$$
\mathbb{P}\left(\sup_{t \in [0,1]^d} |\xi_t| > s\sqrt{(1-\delta)^{1/p} 2d |\log(h_n)|}\right) \geq \sqrt{1-\delta}.
$$

Finally, for $h_n \leq h_0$,

$$
\begin{aligned}
\mathbb{E}\|\xi\|_\infty^p &\geq \mathbb{E}\|\xi\|_\infty^p \mathbf{1}_{\|\xi\|_\infty > s\sqrt{(1-\delta)^{1/p} 2d |\log(h_n)|}} \\
&\geq \left((1-\delta)^{1/p} 2d \|K\|_2^2 |\log(h_n)|\right)^{p/2} \times \mathbb{P}\left(\sup_{t \in [0,1]^d} |\xi_t| > s\sqrt{(1-\delta)^{1/p} 2d |\log(h_n)|}\right) \\
&\geq \left((1-\delta)^{1/p} 2d \|K\|_2^2 |\log(h_n)|\right)^{p/2} \times \sqrt{1-\delta} \\
&\geq (1-\delta) \left(2d \|K\|_2^2 |\log(h_n)|\right)^{p/2}
\end{aligned}
$$

and

$$
\mathbb{E}\|\tilde{f}_{K,h_n} - \mathbb{E}\tilde{f}_{K,h_n}\|_\infty^p = \mathbb{E}\|Z_n\|_\infty^p \geq (1-\delta) \left(\frac{2d\sigma^2 \|K\|_2^2 |\log(h_n)|}{n h_n^d}\right)^{p/2}.
$$

### 3.3 Proof of Theorem 1

Let us first prove that (5) is true. Using Proposition 1 and the inequality

$$\mathbb{E}\|\tilde{f}_{K,h_n} - \mathbb{E}\tilde{f}_{K,h_n}\|_\infty^p \leq 2^p \mathbb{E}\|\tilde{f}_{K,h_n} - f\|_\infty^p$$

(see Lemma 1), we have for $n$ large enough,

$$\frac{\sqrt{|\log(h_n)|}}{\sqrt{nh_n^d}} \leq C \left(\frac{\log n}{n}\right)^{\beta/(2\beta+d)},$$

where $C$ is a constant. It yields

$$d|\log(h_n)| + \log|\log(h_n)| \leq 2\log(C) + \frac{d}{2\beta+d}\log n + \frac{2\beta}{2\beta+d}\log(\log n). \tag{11}$$

If we set

$$u_n = |\log(h_n)| - \frac{\log n}{2\beta+d},$$

then (11) yields that for $n$ large enough,

$$du_n + \log\left(\frac{\log n}{2\beta+d} + u_n\right) \leq 2\log(C) + \frac{2\beta}{2\beta+d}\log(\log n)$$

and

$$u_n + \frac{1}{d}\log\left(1 + \frac{(2\beta+d)u_n}{\log n}\right) \leq -\frac{\log(\log n)}{2\beta+d} + \frac{1}{d}\log((2\beta+d)C^2). \tag{12}$$

This inequality yields that for $n$ large enough, $u_n < 0$. Now, let us assume that (5) is not true. So, there exists an increasing function $\phi$ such that $(h_{\phi(n)})_n$ satisfies

$$h_{\phi(n)}^{-1}\left(\frac{\log(\phi(n))}{\phi(n)}\right)^{1/(2\beta+d)} \xrightarrow{n \to +\infty} +\infty.$$

So,

$$u_{\phi(n)} + \frac{1}{2\beta+d}\log(\log(\phi(n))) \xrightarrow{n \to +\infty} +\infty. \tag{13}$$

Since $u_{\phi(n)} < 0$ for $n$ large enough, $u_{\phi(n)} = O(\log(\log(\phi(n)))) = o(\log(\phi(n)))$. So Inequalities (12) and (13) are contradictory and (5) is true.

Now, we have for any $n \in \mathbb{N}^*$,

$$
\begin{aligned}
h_n^{-p\beta} \|\mathbb{E}\tilde{f}_{K,h_n} - f\|_\infty^p &\leq h_n^{-p\beta} \mathbb{E} \|\tilde{f}_{K,h_n} - f\|_\infty^p \\
&\leq h_n^{-p\beta} C_1 \left(\frac{\log n}{n}\right)^{p\beta/(2\beta+d)} \\
&\leq C_2,
\end{aligned}
$$

where $C_1$ and $C_2$ denote two constants. So, for any $n \in \mathbb{N}^*$,

$$\sup_n \left[ h_n^{-\beta} \|\mathbb{E}\tilde{f}_{K,h_n} - f\|_\infty \right] < \infty.$$

## 3.4 Proof of Theorem 2

The result of this theorem is obtained by doing a balance between the bias and the variance of the estimator $\hat{f}_{n,\beta}$. Denote for $t \in [0,1]^d$ the bias term

$$\forall\, n \in \mathbb{N}^*, \quad b_n(f, t) = \mathbb{E}(\hat{f}_{n,\beta}(t)) - f(t)$$

and its stochastic term

$$\forall\, n \in \mathbb{N}^*, \quad Z_n(t) = \hat{f}_{n,\beta}(t) - \mathbb{E}(\hat{f}_{n,\beta}(t)) = \frac{\sigma}{h_{n,\beta}^d \sqrt{n}} \int_{\mathbb{R}^d} K\left(\frac{t-u}{h_{n,\beta}}\right) dW_u.$$

By (8), the bias satisfies

$$\forall\, n \in \mathbb{N}^*, \quad \|b_n(f, \cdot)\|_\infty \leq C_1 h_{n,\beta}^\beta, \tag{14}$$

where $C_1$ is a constant. Let $A$ such that the support of $K$ is included in $[-A, A]^d$. We consider $n$ large enough to have $Ah_{n,\beta} < 1 - Ah_{n,\beta}$. Here we consider first the case $d = 1$. Since for $d = 1$, $\forall\, t \in [0,1]$, $K_{h_{n,\beta}}(t - \cdot)$ is supported by $[t - Ah_{n,\beta}, t + Ah_{n,\beta}] \subset [-1, 2]$, using Definition (3), we can write $Z_n$ of the form

$$Z_n(t) = Z_n^1(t) + Z_n^2(t) + Z_n^3(t) \tag{15}$$

with

$$Z_n^1(t) = \frac{\sigma}{h_{n,\beta}\sqrt{n}} \int_{[-1,0]} K\left(\frac{t-u}{h_{n,\beta}}\right) dW_u = \frac{\sigma}{h_{n,\beta}\sqrt{n}} \int_{[0,1]} K\left(\frac{u+1-t}{h_{n,\beta}}\right) dW_u,$$

$$Z_n^2(t) = \frac{\sigma}{h_{n,\beta}\sqrt{n}} \int_{[0,1]} K\left(\frac{t-u}{h_{n,\beta}}\right) dW_u,$$

$$Z_n^3(t) = \frac{\sigma}{h_{n,\beta}\sqrt{n}} \int_{[1,2]} K\left(\frac{t-u}{h_{n,\beta}}\right) dW_u = \frac{\sigma}{h_{n,\beta}\sqrt{n}} \int_{[0,1]} K\left(\frac{u-1-t}{h_{n,\beta}}\right) dW_u.$$

Since $K$ satisfies conditions $(A_3)$, we have for any $j \in \{1,2,3\}$ and for all $(s,t) \in [0,1]^2$

$$\mathbb{E}\left[(Z_n^j(t) - Z_n^j(s))^2\right] \leq \frac{C\sigma^2}{nh_{n,\beta}^d} \left\|\frac{t-s}{h_{n,\beta}}\right\|^{2\gamma}, \tag{16}$$

where $C$ is a constant. Moreover we have

$$\mathbb{E}\left[(Z_n^j(t))^2\right] \leq \frac{\sigma^2\|K\|_2^2}{nh_{n,\beta}^d}. \tag{17}$$

Following the same lines of the proof of Lemma 4 in the Appendix of Bertin (2004), using (16) and (17), we have for $r > 0$ and $n$ large enough

$$\mathbb{P}\left[\sup_{t \in [0,1]^d} |Z_n^j(t)| \geq r\sqrt{\frac{2d\sigma^2\|K\|_2^2|\log h_{n,\beta}|}{nh_{n,\beta}^d}}\right] \leq \frac{C_2|\log h_{n,\beta}|^{d/(2\gamma)}}{h_{n,\beta}^d} \exp\{-r^2 d|\log h_{n,\beta}|\}\exp\{C_3 r\},$$

with $C_2$ and $C_3$ positive constants. Then using that for a positive variable $X$ we have $\mathbb{E}[X] = \int_0^{+\infty} \mathbb{P}[X \geq t]dt$, we deduce that for any $j \in \{1,2,3\}$,

$$\forall\, n \in \mathbb{N}^*, \quad \mathbb{E}\left(\|Z_n^j\|_\infty^p\right) \leq C_4 \left(\frac{|\log h_{n,\beta}|}{nh_{n,\beta}^d}\right)^{p/2},$$

with $C_4$ a positive constant. Then,

$$\forall\, n \in \mathbb{N}^*, \quad \mathbb{E}\left(\|Z_n\|_\infty^p\right) \leq C_5 \left(\frac{|\log h_{n,\beta}|}{nh_{n,\beta}^d}\right)^{p/2},$$

with $C_5$ a positive constant.

For $d > 1$, as in (15), $Z_n$ can be decomposed as the sum of $3^d$ terms that satisfy (16) and

(17) and consequently

$$\forall\, n \in \mathbb{N}^*, \quad \mathbb{E}\left(\|Z_n\|_\infty^p\right) \leq C_6 \left(\frac{|\log h_{n,\beta}|}{nh_{n,\beta}^d}\right)^{p/2}, \tag{18}$$

with $C_6$ a positive constant that depends on $d$.

Then with $h_{n,\beta} = C \left(\frac{\log n}{n}\right)^{1/(2\beta+d)}$, (14) and (18) imply (9).

### 3.5 Proof of Theorem 4

First, let us establish the following result concerning the class $\mathcal{I}(\lceil\beta\rceil)$.

**Proposition 3.** *Let* $(\tilde{f}_{K,h_n}) \in \mathcal{I}(\lceil\beta\rceil)$. *If* $f \in \mathbb{L}_\infty^{per}(\mathbb{R}^d)$ *such that*

$$\sup_n \left\{h_n^{-\beta}\|K_{h_n} * f - f\|_\infty\right\} < \infty, \tag{19}$$

*then* $f \in \mathcal{B}_{\infty,\infty}^\beta$.

Remember that for any $n$, $h_n = 2^{-m_n}$. Without loss of generality, we assume that $m_1 = 1$ and $\sup_n (m_{n+1} - m_n) \leq 1$. We have

$$\tilde{f}_{K,h_n}(t) = \frac{1}{h_n^d} \int_{\mathbb{R}^d} K\left(\frac{t-u}{h_n}\right) dY_u$$

with $K \in \mathcal{K}(\lceil\beta\rceil)$ and $(h_n)_n \in \mathcal{H}$. We denote for $t \in [0,1]^d$,

$$b_{h_n}(f,t) = \mathbb{E}\tilde{f}_{K,h_n}(t) - f(t) = K_{h_n} * f(t) - f(t).$$

So,

$$\|b_{h_n}(f,\cdot)\|_\infty \leq C h_n^\beta,$$

where $C$ is a constant. For $k \in \mathbb{N}^*$, we set

$$\begin{aligned}\mathcal{N}_k &= \left\{n: \; 2^{-k-1} < h_n \leq 2^{-k}\right\} \\ &= \left\{n: \; 2^{-k-1} < 2^{-m_n} \leq 2^{-k}\right\} \\ &= \left\{n: \; k \leq m_n < k+1\right\}\end{aligned}$$

For any $k \in \mathbb{N}^*$, $\mathcal{N}_k \neq \emptyset$, and we denote $n_k = \max(\mathcal{N}_k)$ and $h_k^* = h_{n_k}$ so $2^{-k-1} < h_k^* \leq 2^{-k}$. We set

$$u_0 = K_{h_1^*} * K_{h_1^*} * f,$$

and for $k \in \mathbb{N}^*$,

$$u_k = K_{h_{k+1}^*} * K_{h_{k+1}^*} * f - K_{h_k^*} * K_{h_k^*} * f.$$

Using (19), we have

$$\begin{aligned}
\|K_{h_k^*} * K_{h_k^*} * f - f\|_\infty &\leq \|K_{h_k^*} * K_{h_k^*} * f - K_{h_k^*} * f\|_\infty + \|K_{h_k^*} * f - f\|_\infty \\
&= \|K_{h_k^*} * b_{h_k^*}(f, \cdot)\|_\infty + \|b_{h_k^*}(f, \cdot)\|_\infty \\
&\leq \|b_{h_k^*}(f, \cdot)\|_\infty \times \int |K(u)| du + \|b_{h_k^*}(f, \cdot)\|_\infty \\
&\leq C_1 (h_k^*)^\beta \\
&\leq C_1 2^{-k\beta}
\end{aligned}$$

where $C_1$ is a constant. So, we have

$$\lim_{k \to +\infty} \|K_{h_k^*} * K_{h_k^*} * f - f\|_\infty = 0$$

since $\lim_{n \to +\infty} h_n = 0$ and then

$$\sum_{k \geq 0} u_k = f.$$

Furthermore, for $k \geq 1$, since $(h_n)_n \in \mathcal{H}$

$$\begin{aligned}
\|u_k\|_\infty &\leq \|K_{h_k^*} * K_{h_k^*} * f - f\|_\infty + \|K_{h_{k+1}^*} * K_{h_{k+1}^*} * f - f\|_\infty \\
&\leq C_1 (h_k^*)^\beta + C_1 (h_{k+1}^*)^\beta \\
&\leq 2C_1 (h_k^*)^\beta \\
&\leq 2C_1 2^{-k\beta}.
\end{aligned}$$

Let $N = \lceil \beta \rceil$ and $(\alpha_1, \ldots, \alpha_d) \in \mathbb{N}^d$, such that $\sum_{i=1}^d \alpha_i = N$. Using the properties of the convolution operator and condition $(A_5)$ on the kernel, we have for $k \geq 1$,

$$
\begin{aligned}
\left\| \frac{\partial^{\alpha_1}}{\partial t_1^{\alpha_1}} \cdots \frac{\partial^{\alpha_d}}{\partial t_d^{\alpha_d}} u_k \right\|_\infty &= \left\| \frac{\partial^{\alpha_1}}{\partial t_1^{\alpha_1}} \cdots \frac{\partial^{\alpha_d}}{\partial t_d^{\alpha_d}} \left( K_{h^*_{k+1}} * K_{h^*_{k+1}} * f - K_{h^*_k} * K_{h^*_k} * f \right) \right\|_\infty \\
&= \left\| \frac{\partial^{\alpha_1}}{\partial t_1^{\alpha_1}} \cdots \frac{\partial^{\alpha_d}}{\partial t_d^{\alpha_d}} \left( K_{h^*_{k+1}} + K_{h^*_k} \right) * \left( K_{h^*_{k+1}} - K_{h^*_k} \right) * f \right\|_\infty \\
&\leq \left( (h^*_{k+1})^{-N} + (h^*_k)^{-N} \right) \| K_{h^*_{k+1}} * f - K_{h^*_k} * f \|_\infty \int \left| \frac{\partial^{\alpha_1}}{\partial t_1^{\alpha_1}} \cdots \frac{\partial^{\alpha_d}}{\partial t_d^{\alpha_d}} K(t) \right| dt \\
&\leq C_2 (h^*_k)^{-N} \left( \| b_{h^*_{k+1}}(f, \cdot) \|_\infty + \| b_{h^*_k}(f, \cdot) \|_\infty \right) \\
&\leq C_3 (h^*_k)^{-N+\beta} \\
&\leq C_3 2^{(k+1)(N-\beta)},
\end{aligned}
$$

where $C_2$ and $C_3$ are constants. We have used $(h_n)_n \in \mathcal{H}$.

Now, by setting for any $j \geq 1$, $f_j = \sum_{k=0}^{j-1} u_k$, we have

$$\sup_{j \in \mathbb{N}^*} 2^{j\beta} \| f - f_j \|_\infty < \infty$$

and

$$\sup_{j \in \mathbb{N}^*} 2^{-(N-\beta)j} \left\| \frac{\partial^{\alpha_1}}{\partial t_1^{\alpha_1}} \cdots \frac{\partial^{\alpha_d}}{\partial t_d^{\alpha_d}} f_j \right\|_\infty < \infty.$$

These inequalities prove that $f$ belongs to $f \in \mathcal{B}^{\beta}_{\infty,\infty}$ (see Definition 3).    $\square$

Now, let us prove Theorem 4. Theorem 3 and Proposition 3 imply that

$$MS(\hat{f}_\beta, \psi(\beta), p) \subset \mathcal{B}^{\beta}_{\infty,\infty}.$$

Now, let us consider $f \in \Sigma(\beta)$. This implies

$$|f(b) - P_m(f)(b-a, a)| \leq L \| b - a \|^\beta,$$

with $L > 0$, for all $a = (a_1, \ldots, a_d), b = (b_1, \ldots, b_d) \in \mathbb{R}^d$ where $P_m(f)(x, a)$ is the Taylor polynom of order $m$ associated to the function $f$ in the neighborhood of $a$ and $m = \lfloor \beta \rfloor$. For

$t \in [0,1]^d$ and $n$ large enough, we have, since $K$ satisfies $(A_4)$

$$
\begin{aligned}
K_{h_{n,\beta}} * f(t) - f(t) &= \frac{1}{h_{n,\beta}^d} \int_{\mathbb{R}^d} K\left(\frac{t-u}{h_{n,\beta}}\right) (f(u) - f(t)) du \\
&= \frac{1}{h_{n,\beta}^d} \int_{\mathbb{R}^d} K\left(\frac{t-u}{h_{n,\beta}}\right) [f(u) - P_m(f)(u-t,t)] du,
\end{aligned}
$$

where the last line comes from that $K$ satisfies $(A_6)$. We have, for $n$ large enough,

$$
\begin{aligned}
|K_{h_{n,\beta}} * f(t) - f(t)| &\leq L h_{n,\beta}^\beta \frac{1}{h_{n,\beta}^d} \int_{\mathbb{R}^d} K\left(\frac{t-u}{h_{n,\beta}}\right) \frac{\|t-u\|^\beta}{h_{n,\beta}^\beta} du \\
&= L h_{n,\beta}^\beta \int_{\mathbb{R}^d} K(u) \|u\|^\beta du
\end{aligned}
$$

and the last line comes from a change of variables. Using Theorem 3, this implies that

$$\Sigma(\beta) \subset MS(\hat{f}_\beta, \psi(\beta), p).$$

### 3.6 Proof of Theorem 5

Here we emphasize on the fact that since for all $\gamma \in B$, $\gamma$ is not an integer, we have for all $\gamma \in B$ $\Sigma(\gamma) = \mathcal{B}_{\infty,\infty}^\gamma$. Consider the following proposition:

**Proposition 4.** *Let $\beta \in B$. If $f \in \mathcal{B}_{\infty,\infty}^\beta$, then*

$$\sup_n \left\{ (\psi_n(\beta))^{-p} \mathbb{E}\left[\|f - \hat{f}_{n,\hat{\beta}}\|_\infty^p\right] \right\} < \infty. \tag{20}$$

*Proof.* The proof of (20) is the same as the proof of Theorem 2 of Bertin (2005) but for kernels more general and regularities larger than 1. □

Proposition 4 proves that

$$\mathcal{B}_{\infty,\infty}^\beta \subset MS(\hat{f}, \psi(\beta), p).$$

Now, here we prove the inclusion

$$MS(\hat{f}, \psi(\beta), p) \subset \mathcal{B}_{\infty,\infty}^\beta. \tag{21}$$

Let $f \in MS(\hat{f}, \psi(\beta), p)$. This implies that

$$\sup_n \left\{ (\psi_n(\beta))^{-p} \, \mathbb{E}\left[ \|f - \hat{f}_{n,\hat{\beta}}\|_\infty^p \right] \right\} < \infty. \tag{22}$$

Moreover, we have

$$\forall \, \gamma \leq \hat{\beta}, \quad \|\hat{f}_{n,\gamma} - \hat{f}_{n,\hat{\beta}}\|_\infty \leq \eta_n(\gamma). \tag{23}$$

Let us establish some preliminary results.

**Lemma 2.** *If $f \in \mathcal{B}_{\infty,\infty}^\gamma$, then $\mathbb{P}(\hat{\beta} < \gamma) \to 0$, when $n$ goes to $+\infty$.*

*Proof.* The proof is the same as the proof of Lemma 7 of Bertin (2005) and we have for any $\gamma \in B$, if $f \in \mathcal{B}_{\infty,\infty}^\gamma$, then $\mathbb{P}(\hat{\beta} < \gamma) \to 0$ with a polynomial rate when $n$ goes to $+\infty$. $\square$

**Lemma 3.** *Let $f \in MS(\hat{f}, \psi(\beta), p)$. Under the assumptions of Theorem 5, we have $f \in \mathcal{B}_{\infty,\infty}^{\beta_1}$.*

*Proof.* Since $\beta_1 \leq \hat{\beta}$ and $\beta_1 \leq \beta$, using (22) and (23), we have that

$$\begin{aligned}
\mathbb{E}\|\hat{f}_{n,\beta_1} - f\|_\infty^p &\leq 2^{p-1}\left( \mathbb{E}\|\hat{f}_{n,\beta_1} - \hat{f}_{n,\hat{\beta}}\|_\infty^p + \mathbb{E}\|f - \hat{f}_{n,\hat{\beta}}\|_\infty^p \right) \\
&\leq 2^{p-1}\left( (\eta_n(\beta_1))^p + C(\psi_n(\beta))^p \right) \\
&\leq \tilde{C}(\psi_n(\beta_1))^p,
\end{aligned}$$

where $C$ and $\tilde{C}$ are two positive constants. Theorem 4 implies that $f \in \mathcal{B}_{\infty,\infty}^{\beta_1}$. $\square$

Now, here we prove by induction on $\delta \in \{\beta_1, \beta_2, \ldots, \beta\}$ that $f \in \mathcal{B}_{\infty,\infty}^\delta$.
- By Lemma 3, we have that $f \in \mathcal{B}_{\infty,\infty}^{\beta_1}$.
- Now, let $\delta \in \{\beta_2, \beta_3, \ldots, \beta\}$. We assume that $f \in \mathcal{B}_{\infty,\infty}^{\delta^-}$, where

$$\delta^- = \max\{\gamma \in B: \quad \gamma < \delta\}.$$

Now, to obtain that $f \in \mathcal{B}_{\infty,\infty}^\delta$, it is enough to prove that there exists a constant $C$ such that for all $n$ large enough,

$$\|\mathbb{E}\hat{f}_{n,\delta} - f\|_\infty^p \leq C(\psi_n(\delta))^p$$

(see Proposition 3). We have $\forall \, n \in \mathbb{N}^*$,

$$\begin{aligned}
\left(\psi_n(\delta^-)\right)^p &\leq C_1 \mathbb{E}\|\hat{f}_{n,\delta^-} - \mathbb{E}\hat{f}_{n,\delta^-}\|_\infty^p \\
&\leq 2^p C_1 \mathbb{E}\|\hat{f}_{n,\delta^-} - f\|_\infty^p \\
&\leq 2^p C_1 \left(\mathbb{E}\left[\|\hat{f}_{n,\delta^-} - f\|_\infty^p \mathbf{1}_{\hat{\beta}=\delta^-}\right] + \mathbb{E}\left[\|\hat{f}_{n,\delta^-} - f\|_\infty^p \mathbf{1}_{\hat{\beta}\neq\delta^-}\right]\right) \\
&\leq 2^p C_1 \left(\mathbb{E}\|\hat{f}_{n,\hat{\beta}} - f\|_\infty^p + \left(\mathbb{E}\|\hat{f}_{n,\delta^-} - f\|_\infty^{2p}\right)^{1/2} \mathbb{P}^{1/2}(\hat{\beta} \neq \delta^-)\right) \\
&\leq C_2 \left((\psi_n(\beta))^p + (\psi_n(\delta^-))^p \mathbb{P}^{1/2}(\hat{\beta} < \delta^-) + (\psi_n(\delta^-))^p \mathbb{P}^{1/2}(\hat{\beta} > \delta^-)\right) \\
&\leq C_2 \left(\psi_n(\delta^-)\right)^p \mathbb{P}^{1/2}(\hat{\beta} > \delta^-) + o\left((\psi_n(\delta^-))^p\right) \\
&\leq C_2 \left(\psi_n(\delta^-)\right)^p \mathbb{P}^{1/2}(\hat{\beta} \geq \delta) + o\left((\psi_n(\delta^-))^p\right),
\end{aligned}$$

with $C_1$ and $C_2$ two positive constants, where the first two lines are a consequence of Lemma 1 and Proposition 1, the fifth line comes from (22) and the sixth line is a consequence of Lemma 2 and $\delta \leq \beta$. This implies that for $n$ large enough,

$$\mathbb{P}(\hat{\beta} \geq \delta) \geq \frac{1}{4C_2^2}.$$

Now we have for $n$ large enough

$$\begin{aligned}
\|\mathbb{E}\hat{f}_{n,\delta} - f\|_\infty^p &\leq 4C_2^2 \mathbb{P}(\hat{\beta} \geq \delta)\|\mathbb{E}\hat{f}_{n,\delta} - f\|_\infty^p \\
&\leq 4C_2^2 \mathbb{E}\left[\mathbf{1}_{\{\hat{\beta}\geq\delta\}}\|\mathbb{E}\hat{f}_{n,\delta} - f\|_\infty^p\right] \\
&\leq C_3 \mathbb{E}\left[\mathbf{1}_{\{\hat{\beta}\geq\delta\}}\left(\|\mathbb{E}\hat{f}_{n,\delta} - \hat{f}_{n,\delta}\|_\infty^p + \|\hat{f}_{n,\delta} - \hat{f}_{n,\hat{\beta}}\|_\infty^p + \|\hat{f}_{n,\hat{\beta}} - f\|_\infty^p\right)\right],
\end{aligned}$$

where $C_3$ is a positive constant. Using properties (22) and (23) and (3.4), we deduce that

$$\|\mathbb{E}\hat{f}_{n,\delta} - f\|_\infty^p \leq C_4 \left(\psi_n(\delta)\right)^p,$$

with $C_4$ a positive constant.

Then using that the induction is for $\delta \in \{\beta_1, \beta_2, \ldots, \beta\}$, we deduce that $f \in \mathcal{B}^{\beta}_{\infty,\infty}$ and we obtain the inclusion (21).

# References

Autin, F., (2006). Maxisets for density estimation on R (2006), *Math. Methods of Statist.* **15**(2), 123–145.

Autin, F., Picard, D. & Rivoirard, V. (2006). Maxiset approach for Bayesian nonparametric estimation. To appear in *Math. Methods of Statist.*

Bertin, K. (2004). Asymptotically exact minimax estimation in sup-norm for anisotropic Hölder classes, *Bernoulli* **10**(5). 873–888.

Bertin, K. (2005). Sharp adaptive estimation in sup-norm for $d$-dimensional Hölder classes, *Math. Methods Statist.* **14**(3), 267–298.

Cohen, A., DeVore, R., Kerkyacharian, G. & Picard, D. (2001). *Maximal spaces with given rate of convergence for thresholding algorithms*, Appl. Comput. Harmon. Anal. **11**(2), 167–191.

Ibragimov, I.A. & Hasminskii, R. Z. (1981). *Statistical estimation, Asymptotic theory*, Applications of Mathematics 16. Springer-Verlag, New York.

Kerkyacharian, G. & Picard, D. (1993). Density estimation by kernel and wavelets methods: optimality of Besov spaces, *Statist. Probab. Lett.* **18**(4), 327–336.

Kerkyacharian, G. & Picard, D. (2000). Thresholding algorithms, maxisets & well-concentrated bases, *Test* **9**(2), 283–344.

Kerkyacharian, G. & Picard, D. (2002). Minimax or maxisets?, *Bernoulli* **8**(2), 219–253.

Lepskii, O. V. (1992). *On problems of adaptive estimation in white Gaussian noise,* In: Topics in nonparametric estimation, Adv. Soviet Math. 12, Amer. Math. Soc.Providence, RI, 87–106.

Meyer, Y. (1990). *Ondelettes et opérateurs. I.*, Actualités Mathématiques, Hermann, Paris.

Rivoirard, V. (2004). Maxisets for linear procedures, *Statist. Probab. Lett.* **67**(3), 267–275.

Rivoirard, V. (2005). Bayesian modeling of sparse sequences and maxisets for Bayes rules, *Math. Methods Statist.* **14**(3), 346–376.

Rivoirard, V. & Tribouley, K. (2006) Maxisets for integrated quadratic functionals. To appear in *Statistica Sinica*.

Stone, C. J. (1982). Optimal global rates of convergence for nonparametric regression, *Ann. Statist.* **10**(4), 1040–1053.

Tsybakov, A. B. (2004). *Introduction à l'estimation non-paramétrique* Springer-Verlag, Berlin.

Vincent Rivoirard, Equipe Probabilité, Modélisation et Statistique, Laboratoire de Mathématique, CNRS UMR 8628, Université Paris Sud, 91405 Orsay Cedex, France.
E-mail : Vincent.Rivoirard@math.u-psud.fr